\documentstyle[11pt]{amsart}

\newtheorem{thm}{Theorem}[section]
\newtheorem{prop}[thm]{Proposition}
\newtheorem{lem}[thm]{Lemma}
\newtheorem{cor}[thm]{Corollary}
\newtheorem{conjecture}[thm]{Conjecture}

\theoremstyle{definition}
\newtheorem{defn}[thm]{Definition}
\newtheorem{remark}[thm]{Remark}
\newtheorem{example}[thm]{Example}

\theoremstyle{remark}

\numberwithin{equation}{section}

\newcommand{\Wedge}{\bigwedge}

\newcommand{\Span}{\operatorname{Span}}

\newcommand{\Hom}{\operatorname{Hom}}

\newcommand{\bbA}{{\Bbb A}}
\newcommand{\bbG}{{\Bbb G}}

\newcommand{\bbZ}{{\Bbb Z}}

\newcommand{\bbQ}{{\Bbb Q}}

\newcommand{\bbP}{{\Bbb P}}

\newcommand{\cO}{{\cal O}}
\newcommand{\cM}{{\cal M}}

 \newcommand{\brokrarr}{\vphantom{\to}\mathrel{\smash{{-}{\rightarrow}}}}

\newcommand{\lf}{\mathopen}

\let\r=\mathclose

\newcommand{\notdiv}{\mathrel{\not|}}

\newcommand{\sdp}{\mathbin{{>}\!{\triangleleft}}}

\newcommand{\GL}{{\operatorname{GL}}}
\newcommand{\SL}{{\operatorname{SL}}}

\newcommand{\PGL}{{\operatorname{PGL}}}
\newcommand{\PGLn}{{\operatorname{PGL}_n}}

\newcommand{\rank}{\operatorname{rank}}
\newcommand{\Stab}{\operatorname{Stab}}

\newcommand{\lra}{\longrightarrow}
\newcommand{\Pf}{\operatorname{Pf}}

\newcommand{\Mat}{{\operatorname{M}}}
\newcommand{\Mn}{\Mat_n}

\newcommand{\Sym}{{\operatorname{S}}}

\newcommand{\lm}{\operatorname{lm}}

\let\to=\longrightarrow

\begin{document}
\title[Birational invariant] 
{A birational invariant for algebraic group actions}
\author[Z. REICHSTEIN and B. YOUSSIN]
{Z. Reichstein and B. Youssin} 
\address{Department of Mathematics, Oregon State University,
Corvallis, OR 97331}
\thanks{Z. Reichstein was partially supported by NSF grant DMS-9801675}
\email{zinovy@@math.orst.edu}
\address{Department of Mathematics and Computer Science,
University of the Negev, Be'er Sheva', Israel\hfill\break
\hbox{{\rm\it\hskip\parindent Current address\/}}: 
Hashofar 26/3, Ma'ale Adumim, Israel}
\email{youssin@@math.huji.ac.il}

\begin{abstract}
We construct a birational invariant for certain algebraic group actions.
We use this invariant to classify linear representations of 
finite abelian groups up to birational equivalence, thus
answering, in a special case, a question of E. B. Vinberg and
giving a family of counterexamples to a related
conjecture of P. I. Katsylo. We also give a new proof of 
a theorem of M. Lorenz on birational equivalence of quantum tori 
(in a slightly expanded form) by applying our invariant 
in the setting of $\PGLn$-varieties.
\end{abstract}

\maketitle
\tableofcontents

\section{Introduction}
\label{sect1}

Let $G$ be an algebraic group and let $X$ be a smooth projective
$G$-variety (i.e., an algebraic variety with a $G$-action) defined
over an algebraically closed base field of characteristic zero.
It is shown in~\cite{ry1} that for each finite abelian 
subgroup $H$, the presence of an $H$-fixed points is a birational 
invariant of $X$ as a $G$-variety. 
In other words, if $X$ and $Y$ are birationally isomorphic smooth
projective $G$-varieties and $H$ is a finite abelian subgroup of $G$
then $X^H \neq \emptyset$ iff $Y^H \neq \emptyset$.  
(Note that only nontoral finite abelian subgroups $H$
are of interest in this setting. If $H$ lies in a torus $T$ of $G$
then the Borel Fixed Point Theorem says that $X^H$ can never be empty.) 
In~\cite{ry1},~\cite{ry2} and~\cite{tsch} we used 
this fact to study the geometry of $G$-varieties (their essential 
dimensions, splitting degrees, etc.) and properties of 
related algebraic objects (field extensions, division 
algebras, octonion algebras, etc.)

In this paper we will associate (under additional assumptions on
$X$, $G$ and $H$) a more subtle invariant $i(X, x, H)$ to a point
$x \in X^H$; the precise definition 
is given in Section~\ref{sect.def}.  Our main result about
$i(X, x, H)$ is stated below.

Recall that the rank of a finite abelian group $H$ is the minimal number of
generators of $H$ (see Section~\ref{sect2}) and that a $G$-variety $X$ is
called generically free if $\Stab(x) = \{ 1 \}$ for $x$ in general position 
in $X$ (see Section~\ref{sect3}).

\begin{thm} \label{thm1}
Let $G$ be an algebraic group of dimension $d$, $H$ be a finite abelian 
subgroup of $G$ of rank $r$, and $X$, $Y$ be birationally 
isomorphic smooth projective irreducible generically free $G$-varieties 
of dimension $d + r$.
Assume that $\Stab(x)$ is finite for every $x \in X^H$ and
$\Stab(y)$ is finite for every $y \in Y^H$.  Then for every 
$x \in X^H$ there exists a $y \in Y^H$ such that 
$i(Y, y, H) = i(X, x, H)$.
\end{thm}

Informally speaking, the presence of $H$-fixed points $x$ with 
a prescribed value $i(X, x, H)$ (on a suitable model) is a birational 
invariant of $X$ as a $G$-variety. Our proof of Theorem~\ref{thm1}, 
presented in Section~\ref{sect7}, relies on canonical resolution 
of singularities. Note that Theorem~\ref{thm1} 
remains valid even if $X$ and $Y$ are not assumed to be irreducible; 
see Remark~\ref{rem.irred}.

We give two applications of Theorem~\ref{thm1}. 

\subsection*{A birational classification of linear representations}
For our first application recall that by the no-name lemma any 
two generically free linear representations of a given algebraic
group $G$ are stably isomorphic as $G$-varieties; 
see, e.g.,~\cite[1.5.3]{popov}. 
Thus it is natural to try to classify such representations up 
to birational isomorphism.  This problem was proposed 
by E. B. Vinberg~\cite[pp. 494-496]{pv2}; see 
also~\cite[1.5.1]{popov}.  P. I. Katsylo
has subsequently stated the following conjecture: 

\begin{conjecture} 
(\cite{katsylo}; see also~\cite[1.5.10]{popov}) 
\label{conj1}
Let $V$ and $W$ be generically free
linear representations of an algebraic group $G$.
If $\dim(V) = \dim(W)$ then $V$ and $W$ are birationally isomorphic 
as $G$-varieties.
\end{conjecture}

In this paper we will establish the following birational classification of   
faithful linear representations of a diagonalizable group. Recall that
every diagonalizable group $G$ can be uniquely written in the form 
\begin{eqnarray} \label{e.fin-abel}
G = \bbG_m(n_1) \times \dots \times \bbG_m(n_r) & \text{such that} & 
\nonumber \\
\bbG_m(n_1) \subset \dots \subset \bbG_m(n_r) & \text{and each $n_i = 0$
or $\geq 2$}  \, , & 
\end{eqnarray}
where $\bbG_m = k^* = \bbG_m(0)$ denotes the 1-dimensional torus,
$\bbG_m(n) \simeq \bbZ/n\bbZ$ is the $n$-torsion subgroup of $\bbG_m$, and 
$\bbG(a) \subset \bbG(b)$ iff either $a$ divides $b$ or $a = b = 0$.
Recall that a representation of a diagonalizable group is
faithful iff it is generically free and that any such representation
has dimension $\geq r$. 

\begin{thm} \label{thm6.2} 
Let $G$ be a diagonalizable group, as in~\eqref{e.fin-abel}.

\smallskip
(a) If $d \geq r + 1$ then any two faithful $d$-dimensional linear 
representations of $H$ are birationally equivalent.

\smallskip
(b) If $n_1 = 0$ or $2$ then any two faithful
$r$-dimensional linear representations 
of $H$ are birationally equivalent.

\smallskip
(c) If $n_1 \geq 3$ then $H$ has exactly $\phi(n_1)/2$ birational equivalence 
classes of faithful $r$-dimensional representations. Here $\phi$ denotes the
Euler $\phi$-function. 

\smallskip
\noindent
In particular, Conjecture~\ref{conj1} fails for $G$ if and only if
$n_1 = 5$ or $\geq 7$.
\end{thm}

The birational equivalence classes of faithful 
linear representations of $G$ are explicitly described in 
Theorem~\ref{thm6.1}. Later in Section~\ref{sect.bir-class}
we will show that Conjecture~\ref{conj1} also 
fails for a class of nonabelian 
groups $G$. On the other hand, we remark that
P. I. Katsylo~\cite{katsylo} proved Conjecture~\ref{conj1} for
$G = \SL_2$, $G = \PGL_2$, $G = S_n$ ($n \leq 4$) and that
many interesting cases remain open, including
$G = \Sym_n$ ($n \geq 5$) and $G$ = arbitrary 
connected semisimple group. 

\subsection*{Birational equivalence of quantum tori}

Our second application is based on the fact that a generically free
$\PGLn$-varieties $X$ with $k(X)^{\PGLn} = K$ are in 1-1 correspondence 
with division algebras $D$ with center $K$; see, e.g.,~\cite[X.5]{serre}
or~\cite[Section 3]{ry2}. Thus 
Theorem~\ref{thm1} (with $G = \PGLn$) will sometimes allow us 
to prove that certain division algebras are not isomorphic 
to each other over $k$.  

Let $\omega_1, \dots, \omega_r$ be roots of unity and let
$R(\omega_1, \dots, \omega_r)$ be the $k$-algebra
$k\{x_1^{\pm 1}, \dots, x_{2r}^{\pm 1}\}$, 
where $x_{2i-1} x_{2i} = \omega_i x_{2i} x_{2i-1}$ 
for $i =1, \dots, r$ and all other pairs of variables commute.
Denote the algebra of quotients of $R(\omega_1, \dots, 
\omega_r)$ by $Q(\omega_1, \dots, \omega_r)$.  
Note that $Q(\omega_1, \dots, \omega_r)$ is obtained from  
$R(\omega_1, \dots, \omega_r)$ by adjoining the inverses of all
central elements and that $Q(\omega_1, \dots, \omega_r)$ is 
a finite-dimensional division algebra (in fact, it is a tensor product 
of symbol algebras).  

M. Lorenz~\cite[Proposition 1.3]{lorenz} showed that
$Q(\omega)$ and $Q(\omega')$ are isomorphic as $k$-algebras if and only if 
$\omega' = \omega^{\pm 1}$. In Section~\ref{sect.pgl_n}
we will give a geometric proof of the following variant of this result.

\begin{thm} \label{thm3}
Suppose $\omega_i$ is a primitive $n_i$th root of unity, $n_i$ 
divides $n_{i+1}$ for $i = 1, \dots, r-1$, $n_1 \geq 2$, and
$(m_i, n_i) = 1$.
Then $Q(\omega_1, \dots, \omega_r)$ and $Q(\omega_1^{m_1}, \dots,
\omega_r^{m_r})$ are isomorphic as $k$-algebras
if and only if $m_1 ... m_r \equiv \pm1 \pmod{n_1}$.
\end{thm}

The algebras $Q(\omega_1, \dots, \omega_r)$ 
and $Q(\omega_1^{m_1}, \dots, \omega_r^{m_r})$ have a common 
center $K = k(x_1^{n_1}, \dots, x_{2r}^{n_{2r}})$. Note that 
these algebras may be $k$-isomorphic but not $K$-isomorphic,
i.e., not Brauer equivalent.  More precisely, 
$Q(\omega_1, \dots, \omega_r)$ and $Q(\omega_1^{m_1}, \dots,
\omega_r^{m_r})$ are Brauer equivalent iff $m_i \equiv 1 \pmod{n_i}$ 
for every $i = 1, \dots, r$.

It is natural to think of $R(\omega_1, \dots, \omega_r)$ and
$Q(\omega_1, \dots, \omega_r)$ as, respectively, the ``coordinate ring"
and the ``function field" of a quantum torus. Using this terminology,
one may view Theorem~\ref{thm3} as a result about birational 
isomorphism classes of quantum tori. 

Finally we remark that M. Lorenz has communicated to us a proof
of Theorem~\ref{thm3} based on the techniques 
of~\cite{lorenz} and~\cite{lorenz2}. His argument works 
in arbitrary characteristic. 

\subsection*{Acknowledgements}
We are grateful to P. I. Katsylo, M. Lorenz and V. L. Popov
fo helpful communications related to the subject matter of this paper.

The second author warmly thanks Institute of Mathematics of Hebrew
University for its hospitality during 1999/2000.

\section{Linear algebra in abelian groups}
\label{sect2}

Recall that any finitely generated abelian group $(A, +)$ 
can be written in the form
\begin{equation} \label{e.fg-abel}
\begin{gathered}
A \simeq \bbZ/n_1 \bbZ \times \dots \times \bbZ/n_r \bbZ \; ,
\text{ where each $n_i = 0$ or $\geq 2$ and}   \\ 
\text{for every $i = 1, \dots, r-1$ either $n_i = n_{i+1} = 0$
or $n_i$ divides $n_{i+1}$.} 
\end{gathered}
\end{equation}
Here $r$ and $n_1, \dots, n_r$ are uniquely
determined by the isomorphism type of $A$.
We shall refer to the integer $r$ as the {\em the rank} of $A$;
equivalently, the rank of $A$ equals the minimal possible number of 
generators of $A$. 

Recall that if $B$ is an abelian group then the dual group $B^*$ is 
defined as $\Hom(B, \bbQ/\bbZ)$; we will often identify 
$\bbQ/\bbZ$ with the multiplicative group of roots of unity in $k$. 
The finitely generated group $A$ of~\eqref{e.fg-abel}
and the diagonalizable group
$G$ of~\eqref{e.fin-abel} are dual to each other. 
The rank of a diagonalizable group $G$ is defined as the rank
of the finitely generated group $G^*$ (in particular, the group
$G$ of~\eqref{e.fin-abel} has rank $r$).
Note that this is consistent with the way we defined rank for a finitely 
generated group: indeed, if $A$ is both diagonalizable and finitely 
generated, i.e., is finite abelian, then $A$ and $A^*$ are isomorphic, so that
their ranks coincide. 

\subsection*{Skew-symmetric powers}

We will write $\bigwedge^d(A)$ for the $d$th skew-symmetric power
of $A$, viewed as a $\bbZ$-module. 

The proof of the following lemma is elementary; we leave it as 
an exercise for the reader.

\begin{lem} \label{lem.ab1} Let $A$ be a finitely generated 
abelian group as in~\eqref{e.fg-abel}. Then

\smallskip
(a) $\bigwedge^r(A) \simeq \bbZ/n_1 \bbZ$.

\smallskip
(b) $\bigwedge^d(A) = (0)$, if $d \geq r + 1$.
\qed
\end{lem}

\begin{defn} \label{def.sympl} Let $A$ be a finite abelian group.
Let $\omega \colon A \times A \lra \bbQ/\bbZ$ be a $\bbZ$-bilinear form.
As usual, we shall say that

\smallskip
(a) $\omega$ is {\em alternating} if $\omega(a, a) = 0$
for every $a \in A$, 

\smallskip
(b) $\omega$ is {\em non-degenerate} if 
$\omega(a, \cdot)$ is not identically zero for 
any $a \in A - \{ 0 \}$, 

\smallskip
(c) $\omega$ is {\em symplectic} if it is both alternating 
and non-degenerate.
\end{defn}

\begin{lem} \label{lem.sf2} Let $A$ be a finite abelian group of
rank $r$, $\omega$ be a symplectic form of $A$, and $\psi$ be 
an $\omega$-preserving automorphism $A \lra A$.  Then 

\smallskip
(a) $\wedge^{2r} \psi$ is the trivial automorphism 
of $\bigwedge^{2r}(A)$ and

\smallskip
(b) $\wedge^{2r} \psi^*$ is the trivial automorphism 
of $\bigwedge^{2r}(A^*)$.
\end{lem}

\begin{pf} (a) It is well-known that $A$ can be written in the form
$A = A_0 \oplus A_0^*$ such that 
\begin{equation} \label{e.omega}
\omega((a, a^*), (b, b^*)) = a^*(b) - b^*(a)
\end{equation}
for any $a, b \in A_0$ and $a^*, b^* \in A_0^*$; 
see, e.g.,~\cite[Theorem 4.1]{ta}. 
Write $A_0$ as $\bbZ/n_1\bbZ \times \dots \times \bbZ/n_r\bbZ$, 
where $n_i$ divides $n_{i+1}$ for
$i = 1, \dots, r-1$ and $n_1 \geq 2$.
Let $e_i \in A_0$ be a generator of the of the $i$th factor, 
and let $f_i \in A_0^*$ be given by 
$f_i \colon A_0 \lra \bbZ/n_i\bbZ \lra \bbQ/\bbZ$, where the first map 
the projection to the $i$th factor, and the second map takes $e_i$ 
to $1/n_i$.  Then every $a \in A$ can be written in the form 
$a = \sum_{i=1}^n (\alpha_{2i-1} e_i + \alpha_{2i} f_i)$, where
$\alpha_{2i-1}, \alpha_{2i} \in \bbZ$,
and~\eqref{e.omega} translates into 
\[ \omega[\sum_{i=1}^n (\alpha_{2i-1} e_i + \alpha_{2i} f_i),
\sum_{i=1}^n (\beta_{2i-1} e_i + \beta_{2i} f_i)] 
= \sum_{i=1}^n \dfrac{1}{n_i}
(\alpha_{2i} \beta_{2i-1} - \alpha_{2i-1} \beta_{2i}) \, . \]
Suppose 
\[ \begin{array}{l}
 \psi(e_1) = c_{11} e_1 + c_{12} f_1 + \dots + c_{1 , 2r-1} e_r  + 
c_{1 , 2r-1} f_r \, , \\
\psi(f_1) = c_{21} e_1 + c_{22} f_1 + \dots + c_{2 , 2r-1} e_r  + 
c_{2 , 2r} f_r \, ,  \\
\vdots  \\
\psi(e_r) = c_{2r-1 , 1} e_1 + c_{2r-1 , 2} f_1 + \dots + 
c_{2r-1 , 2r-1} e_r  + c_{2r-1 , 2r} f_r \, , \\
\psi(f_r) = c_{2r , 1} e_1 + c_{2r , 2} f_1 + \dots + 
c_{2r , 2r-1} e_r  + c_{2r , 2r} f_r \, , \end{array} \]
where $C = (c_{ij})_{i, j=1, \dots, 2r} \in \Mn(\bbZ)$.
Since $\lambda = e_1 \wedge f_1 \wedge \dots \wedge e_r \wedge f_r$ 
generates $\bigwedge^{2r}(A^*) \simeq \bbZ/n_1 \bbZ$ 
and $[\wedge^{2r}\psi](\lambda) = \det(C) \lambda$, 
it is sufficient to show that $\det(C) = 1 \pmod{n_1}$. 

The condition that $\psi$ preserves $\omega$ 
translates into $C J C^{t} = J \pmod{1}$, 
where $C^t$ is the transpose of $C$ and
\[ J = \left( \begin{array}{ccccc} 0 & 1/n_1 & \dots & 0 & 0 \\ 
                                   -1/n_1 & 0 & \dots & 0 & 0 \\ 
                               \vdots &   &  \vdots &   & \vdots \\       
                                    0 & 0 &        &  0 & 1/n_r  \\ 
                                    0 & 0 &        &  -1/n_r & 0  \end{array} 
\right)  \, .\]
In other words, 
\begin{equation} \label{e.CJC^t}
CJ C^{t} = J + N \, ,
\end{equation}
where $N$ is a skew-symmetric integral matrix. 
We shall deduce the desired equality,
$\det(C) = 1 \pmod{n_1}$, by computing the Pfaffian on both sides
of~\eqref{e.CJC^t}.
On the one hand 
\begin{equation} \label{e.pf-det}
\Pf(CJC^t) = \det(C) \Pf(J)  = (-1)^r \det(C) \dfrac{1}{n_1 \dots n_r} \, ;
\end{equation}
see, e.g.,~\cite[XIV, 10, Theorem 7]{lang}. On the other hand, suppose
$X = (x_{ij})$, where $x_{ji} = - x_{ij}$ for $1 \le i, j \le 2r$. Then 
$\Pf(X) \in \bbZ[x_{ij} \, | \, 1 \leq i < j \leq 2r]$ has degree 1 in every
$x_{ij}$, where $i < j$. (Indeed, $\det(X)$ has degree 2 
in every $x_{ij}$, and $\Pf(X)^2 = \det(X)$.) Consequently, 
\begin{equation} \label{e.pf(J+N)}
 \Pf(J + N) = \Pf(J) + \dfrac{z}{n_2 \dots n_r} 
= (-1)^r \dfrac{1}{n_1 \dots n_r} + \dfrac{z}{n_2 \dots n_r} \, , 
\end{equation}
where $z$ is an integer. (Here we used the fact that $n_i$ 
divides $n_{i+1}$ for every $i = 1, \dots, r-1$.) 
Putting~\eqref{e.CJC^t}, \eqref{e.pf-det} and~\eqref{e.pf(J+N)} together,
we see that $\det(C) = 1 + (-1)^r n_1 z$, i.e., $\det(C) = 1 \pmod{n_1}$,
as claimed.

\smallskip
(b) Let $i \colon A \lra A^*$ be the isomorphism 
$a \mapsto i_a$, where $i_a(b) = \omega(a, b)$. It is easy to see that 
the automorphism $\psi^*: A^* \lra A^*$ preserves the symplectic 
form $\omega^*$ on $A^*$ given by $\omega^*(a^*, b^*) =
\omega(i^{-1}a^*, i^{-1}b^*)$. Applying part (a) to $\psi^*$, 
we obtain the desired result.
\end{pf}

\subsection*{Elementary operations}

Let $A$ be an abelian group. We will say 
that two $d$-tuples $(a_1, \dots, a_d)$ and
$(b_1, \dots, b_d) \in A^d$ are related by 
an {\em elementary operation} if 
$(b_1, \dots, b_d) = 
(a_1, \dots, a_{i-1}, a_i + ma_j, a_{i+1}, \dots, a_d)$
for some $i \neq j$ and $m \in \bbZ$. 

\begin{lem} \label{lem.equiv}
Let $a_1, \dots, a_d$ and $b_1, \dots, b_d$ be two sets of generators
for an abelian group $A$. Then $a = (a_1, \dots, a_d)$ and
$(b_1, \dots, b_d)$ are related by 
a sequence of elementary operations if and only if 
$a_1 \wedge \dots \wedge a_d = b_1 \wedge \dots \wedge b_d$
in $\bigwedge^d(A)$.
\end{lem}

\begin{pf} It is clear that if $(a_1, \dots, a_d)$ and $(b_1, \dots, b_d)$
are related by a sequence of elementary operations then
$a_1 \wedge \dots \wedge a_d = b_1 \wedge \dots \wedge b_d$.
We will prove the converse by induction on $d$.

If $d = 1$ there is nothing to prove,
since $\bigwedge^1(A) = A$. For the induction step,
assume $d \geq 2$ and $A \simeq (\bbZ/n_1\bbZ) 
\times \dots \times (\bbZ/n_r\bbZ)$ is as in~\eqref{e.fg-abel} 
Here $r = \rank(A) \leq d$, since we are assuming $A$ is generated
by $d$ elements.

A $d$-tuple of generators $(a_1, \dots, a_d)$ of $A$ can now be represented
by a $d \times r$-matrix $a = (a_{ij})$ whose $i$th row is $a_i$.
Elements of the $j$th column of this matrix lie in $\bbZ/n_j \bbZ$.
Elementary operation on such matrices amount to adding an integer 
multiple of the $j$th row to the $i$th row or some $i \neq j$.

Elementary operations allow us to perform the Euclidean algorithm 
in the last column of
$(a_{ij})$. Since $a_{1r}, \dots, a_{dr}$ generate $\bbZ/n_r \bbZ$, after
a sequence of elementary operations, we may assume that
\[ (a_{ij}) =  \left( \begin{array}{cccc}  
                                           &      &   & 0 \\ 
                                           & X    &   & \vdots \\ 
                                           &      &   & 0 \\ 
                                         * & \dots& * &  1  
                        \end{array} \right) \, , \]
where $X$ is a $(d-1) \times (r-1)$-matrix. 
Since the rows of $(a_{ij})$
generate $A$, the rows of $X$ generate $\overline{A} = (\bbZ/n_1\bbZ) 
\times \dots \times (\bbZ/n_{r-1}\bbZ)$. Thus
after performing additional elementary operations, we may assume
\[ (a_{ij}) = \left( \begin{array}{cccc}  
                                           &      &   & 0 \\ 
                                           & X    &   & \vdots \\ 
                                           &      &   & 0 \\ 
                                         0 & \dots& 0 &  1  
                        \end{array} \right) \, ,  \; \, 
\text{and similarly} \; \,
 (b_{ij}) = \left( \begin{array}{cccc}  
                                           &      &   & 0 \\ 
                                           & Y    &   & \vdots \\ 
                                           &      &   & 0 \\ 
                                         0 & \dots& 0 &  1  
                        \end{array} \right) \, . \]
In other words, we may assume $(a_1, \dots, a_{d-1})$ and 
$(b_1, \dots, b_{d-1})$ are $(d-1)$-tuples of generators in $\overline{A}$
and $a_d = b_d$ is the generator $1\in\bbZ/n_r \bbZ$. 

We claim that 
\begin{equation} \label{e.d-1}
\textstyle
a_1 \wedge \dots \wedge a_{d-1} = b_1 \wedge \dots 
\wedge b_{d-1} \; \, \text{in $\bigwedge^{(d-1)}(\overline{A})$.}
\end{equation}
Indeed, if $d \geq r+1$, this is obvious, 
since $\bigwedge^{d-1}(\overline{A}) = (0)$; see Lemma~\ref{lem.ab1}(b).
If $d = r$ then the isomorphism $\bigwedge^{d}(A) \simeq \bbZ/n_1\bbZ$ 
identifies $a_1 \wedge \dots \wedge a_d$ with $\det(a_{ij}) \pmod{n_1}$,
and the isomorphism
$\bigwedge^{(d-1)}(\overline{A}) \simeq \bbZ/n_1\bbZ$ identifies
$a_1 \wedge \dots \wedge a_{d-1}$ with $\det(X) \pmod{n_1}$.
Since $a_1 \wedge \dots \wedge a_d = b_1 \wedge \dots \wedge b_d$, we 
know that $\det(a_{ij}) = \det(b_{ij}) \pmod{n_1}$; hence, 
$\det(X) = \det(Y) \pmod{n_1}$, and~\eqref{e.d-1} follows.

Now by the induction assumption 
$(a_1, \dots, a_{d-1})$ and $(b_1, \dots, b_{d-1})$ are related
by a sequence of elementary transformations. Since $a_d = b_d$, so are
$(a_1, \dots, a_d)$ and $(b_1, \dots, b_d)$, as desired.
\end{pf}

\begin{cor} \label{cor.equiv} 
Let $a_1, \dots, a_d$ and $b_1, \dots, b_d$ be two
sets of generators for an abelian group $A$. Then the following conditions
are equivalent.

\smallskip
(a) There exists a matrix $N = (n_{ij}) \in \GL_d(\bbZ)$ such that
$b_i = n_{i1} a_1 + \dots + n_{id}a_d$ for $i = 1, \dots, d$.

\smallskip
(b) $a_1 \wedge \dots \wedge a_d = \pm b_1 \wedge \dots \wedge b_d$
in $\bigwedge^d(A)$.
\end{cor}

\begin{pf} The implication (a) $\Longrightarrow$ (b) is obvious.
To prove the converse, note that we may assume without loss of generality
that $a_1 \wedge \dots \wedge a_d = b_1 \wedge \dots \wedge b_d$; indeed,
if $a_1 \wedge \dots \wedge a_d = - b_1 \wedge \dots \wedge b_d$ then
we can simply replace $(a_1, a_2, \dots, a_d)$
by $(-a_1, a_2, \dots, a_d)$. 
Now Lemma~\ref{lem.equiv} says, $(b_1, \dots, b_d)$ is obtained
from $(a_1, \dots, a_d)$ by a sequence of elementary operations.
Each elementary operation relates the two sets of generators as in (a),
with $N$ = elementary matrix $\in \SL(\bbZ)$. Multiplying these matrices 
we obtain the desired conclusion.
\end{pf} 

\section{$H$-slices}
\label{sect3}

\begin{defn} \label{def.sl1} Let $G$ be an algebraic group and
$X$ be a $G$-variety.  We will call a locally closed subvariety $V$ of $X$ 
a {\em slice} at $x \in V$ if $X$ and $V$ are smooth at $x$ and
$T_x(X) = T_x(Gx) \oplus T_x(V)$. If, moreover, $V$ is invariant
under the action of a subgroup $H$ of $\Stab(x)$, we will 
refer to $V$ as an $H$-slice.
\end{defn}

\begin{remark} \label{rem.sl1.5}
Note that since $V$ is smooth at $x$, we may always replace $V$ 
by its (unique) irreducible component passing through $x$ and thus
assume that $V$ is irreducible.
\end{remark}

\begin{example} \label{ex.induced}
Let $G$ be an algebraic group,
$H$ be an algebraic subgroup of $G$ and $V$ be an $H$-variety. 
Recall that the homogeneous fiber product
$X = G *_H X$ is defined as the geometric quotient
$X = (G \times V)/H$, where $H$ acts on $G \times V$ by $h(g, v) = 
(gh^{-1}, hv)$ (This geometric quotient exists under 
mild assumptions on $V$; in particular, it exists whenever 
$V$ is quasi-projective; see~\cite[Theorem 4.19]{pv}.)
Note that $X = G *_H V$ is naturally a $G$-variety, where $G$
acts by left multiplication on the first factor; the details 
of this construction can be found in~\cite[Section 4.8]{pv}.

The points of $X$ are in 1-1 correspondence with
$H$-orbits in $G \times V$; we shall denote the
point $x \in X$ corresponding to the $H$-orbit of $(g, v)$ 
in $G \times V$ by $x = [g, v]$.
Note that there is an $H$-equivariant map $V \lra X$ given by
$v \mapsto [1, v]$; we shall denote the image of this map by $\tilde{V}$.
With these notations, $\tilde{V}$ 
is an $H$-slice for $X$ at $x = [1, v]$ for every smooth point
$v$ of $V$; see~\cite[Proposition 4.22]{pv}. 
\qed
\end{example}

The following three results, Lemma~\ref{lem.sl3}, Lemma~\ref{lem.sl2} and
Proposition~\ref{prop.faithful}, are immediate consequences 
of the Luna Slice Theorem (see~\cite{luna} or~\cite[Section 6]{pv}), 
if we assume that
$X$ is an affine variety and the orbit of $x$ is closed in $X$.
The arguments we present below do not require these assumptions.

\begin{lem} \label{lem.sl3} Let $G$ be an algebraic group,
$X$ be an irreducible $G$-variety, and
$V$ be a slice at $x \in X$.  Then $GV$ is dense in $X$.
\end{lem}

\begin{pf} Consider the map $\phi \colon G \times V \lra X$, given by
$\phi(g, v) = gv$.
The differential $d\phi_{(1, x)} \colon T_{1}(G) \times
T_x(V) \lra T_x(X)$ is onto, since its image contains both $T_x(Gx)$
and $T_x(V)$.
Consequently, $d\phi$ is onto at a general point of 
$G \times V$. Thus $\phi$ is dominant, and the lemma follows.
\end{pf} 

\begin{lem} \label{lem.sl2} Let $G$ be an algebraic group, 
$H$ be a subgroup, $X$ be a $G$-variety, 
and $x$ is a smooth $H$-fixed point of $X$.
If $H$ is reductive then $X$ has an $H$-slice at $x$.
\end{lem}

\begin{pf} 
Let $\cM_{x, X}$ be the maximal ideal of the local ring of $X$ at $x$.
Consider the natural $H$-equivariant linear maps 
$\cM_{x, X} \lra T_x(X)^* \lra T_x(Gx)^*$.
Since $H$ is reductive these maps have $H$-equivariant splittings
as maps of $k$-vector spaces. Thus we can choose a local  
coordinate system $u_1, \dots, u_n$ in $\cM_{x, X}$ such that
$\Span_k(u_1, \dots, u_d)$ 
is an $H$-invariant $k$-vector subspace of $\cM_{x, X}$ 
and $u_1, \dots, u_d$ (restricted to $Gx$) form a local coordinate system 
in $\cO_{x, Gx}$. (Here $n = \dim(X)$ and $d = \dim(Gx)$.)

Note that $u_1, \dots, u_n$ are regular functions in some $H$-invariant
open neighborhood of $x$.  In this neighborhood a slice $V$ with desired 
properties is given by $u_1 = \dots = u_d = 0$.
\end{pf}

Recall that a $G$-variety $X$ is called generically free 
if $\Stab(x) = \{ 1 \}$ for $x$ in general position in $X$.

\begin{prop} \label{prop.faithful}
Let $G$ be an algebraic group, $H$ be a reductive subgroup,
$X$ be a generically free
$G$-variety and $x$ be a smooth $H$-fixed point of $X$.
Then $H$ acts faithfully on $T_x(X)/T_{x}(Gx)$.
\end{prop}

\begin{pf}  
Let $X_0$ be the unique component of $X$ passing through $x$ 
and $G_0$ be the subgroup of all elements of $G$ that preserve $X_0$.
Note that $G_0$ has finite index in $G$ and $H \subset G_0$. 
After replacing $X$ by $X_0$ and $G$ by $G_0$, we may assume $X$ 
is irreducible.

We now argue by contradiction. Assume the kernel $K$ of the $H$-action on
$T_x(X)/T_x(Gx)$ is non-trivial. Note that $K$ is a normal subgroup
of $H$. Since $H$ is reductive, $K$ is not unipotent. Hence, we can find
a non-identity element $g \in K$ of finite order.

By Lemma~\ref{lem.sl2} $X$ has an $H$-slice $V$ at $x$.  Since 
$T_x(V) \simeq T_x(X)/T_x(Gx)$ as $H$-representations, $g$ acts trivially
on $T_x(V)$. This implies that $g$ acts trivially on $V$; 
see, e.g.,~\cite[Lemma 4.2]{ry1}.
On the other hand, by Lemma~\ref{lem.sl3} $GV$ is dense in $X$; 
consequently, for every $x \in X$ in general position 
$\Stab(x)$ contains a conjugate of $g$. Thus 
the $G$-action on $X$ is not generically free, 
contradicting our assumption.
\end{pf}

\section{Definition and first properties of $i(X, x, H)$}
\label{sect.def}

Throughout this section we shall make the following assumptions: 
\[ \begin{array}{lcl}
 G          & & \mbox{algebraic group} \\
 H          & & \mbox{finite abelian subgroup of $G$ of rank $r$} \\
 X         & & \mbox{$G$-variety of dimension $\dim(G) + r$} \\
 x       & & \mbox{$H$-fixed point of $X$ whose stabilizer is finite} \\
\end{array} \] 

\begin{defn} \label{def-i(x)} 
The $H$-representation on $T_x(X)/T_x(Gx)$ decomposes
as a direct sum of $r$ characters $\chi_1, \dots, \chi_r \in H^*$. We
define 
\[ \textstyle
i(X, x, H) = \chi_1 \wedge \dots \wedge \chi_r \in \bigwedge^r(H^*) 
\, . \]
Since the collection of characters $\chi_1, \dots, \chi_n$ is 
well-defined, up to reordering, $i(X, x, H)$ is well-defined 
in $\Wedge^r(H^*)$, 
up to multiplication by $-1$. Thus, properly speaking, 
$i(X, x, H)$ should be
viewed as an element of the factor set $\Wedge^r(H^*)/\sim$, 
where $w_1 \sim w_2$ iff $w_1 = \pm w_2$.  By abuse of notation we will
sometimes write $i(X, x, H) \in \Wedge^r(H^*)$; in such cases it should 
be understood that $i(X, x, H)$ is only defined up to sign.
\end{defn}

\begin{remark} \label{rem2.25} It is clear from the definition that if
$V$ is an $H$-slice for $X$ at $x$ then $i(X, x, H) = i(V, x, H)$. 
In particular, in the setting of Example~\ref{ex.induced}, if
$V$ is an $r$-dimensional $H$-variety, 
$X = G *_H V$, and $v$ is a smooth $H$-fixed point
of $V$ then
$i(X, [1, v], H) = i(\tilde{V}, [1, v], H) = i(V, v, H)$.
\end{remark}

\begin{remark} \label{rem2.4}
Let $g$ be an element of the normalizer $N_{G}(H)$ and let 
$\phi_g$ be the automorphism of $H$ sending $h$ to $ghg^{-1}$.
Then it is easy to see that 
$i(X, gx, H) = (\wedge^{r}\phi_g^*) \Bigl( i(X, x, H) \Bigr)$,
where $\wedge^r \phi_g^*$ is the automorphism of $\bigwedge^r(H^*)$
induced by $\phi_g$. 
\end{remark}

\begin{example} \label{ex.p^r} 
Let $G = H$ be a finite abelian group of rank $r$, 
$\chi_1, \dots, \chi_r$ be a generating set for $H^*$, and
$V = \bbA^r$ be a faithful $r$-dimensional linear representation 
of $H$, given by 
\[h \colon (v_1, \dots, v_r) \lra (\chi_1(h)v_1, \dots, \chi_r(h)v_r) \, .\]
Then the origin $0_V$ is the only $H$-fixed point of $V$, and
Definition~\ref{def-i(x)} immediately implies
$i(V, 0_V, H) = \chi_1 \wedge \dots \wedge \chi_r$. 
The extended $H$-action on $\overline{V}= \bbP^r$, given by
\[ h(v_0: v_1: \dots: v_d) = (v_0 : \chi_1(h) v_1: \dots: 
\chi_d(h) v_d) \, , \]
has exactly $r+1$ fixed points:
\[ x_0 = (1:0: \dots: 0), \dots, x_r = (0: \dots : 0: 1) \, . \]
Note that $x_0 = 0_V$.  We claim that, up to sign, 
\[ i(\overline{V}, x_j, H) = i(V, x_0, H) = 
\chi_1 \wedge \dots \wedge \chi_r \] 
for $j = 1, \dots, r$. To prove this claim, say for $j = 1$,
note that $v_0/v_1, v_2/v_1, \dots, v_r/v_1$ form an affine coordinate
system on $\overline{V}$ near $x_1$.  The $H$-action
is diagonal in these coordinates, and the representation
of $H$ on $T_{x_1}(\overline{V})$ is the direct sum of the characters 
$\chi_1^{-1}, \chi_2\chi_1^{-1}, \dots, \chi_r \chi_1^{-1}$. 
Consequently,
$i(\overline{V}, x_1, H) = \pm \chi_1^{-1} \wedge \chi_2 \chi_1^{-1} 
\wedge \dots \wedge 
\chi_r \chi_1^{-1} = \pm \chi_1 \wedge \dots \wedge \chi_r$, 
as claimed. 
\end{example}

\begin{lem} \label{lem.stand}
Suppose $X$ is a generically free $G$-variety. Then 
$i(X, x, H)$ generates $\Wedge^r(H^*)$ as a group.
\end{lem}

Note that the statement of the lemma makes sense, even though
$i(X, x, H)$ is only defined up to sign: if $a$ generates $\Wedge^r(H^*)$ 
then so does $-a$.

\begin{pf} By Proposition~\ref{prop.faithful}
$H$ acts faithfully on $T_x(X)/T_x(Gx)$. Hence, the characters
$\chi_1, \dots, \chi_r$ introduced in Definition~\ref{def-i(x)} 
generate $H^*$ as an abelian group, and the lemma follows.
\end{pf}

\section{$i(X, x, H)$ and birational morphisms}
\label{sect6}

The purpose of this section is to prove the following: 

\begin{thm} \label{thm1-1}
Let $G$ be an algebraic group of dimension $d$, 
$H$ be a finite abelian subgroup of $G$ of rank $r$,
$f \colon X \lra Y$ be birational morphism of irreducible 
generically free $G$-varieties of dimension $d+r$, 
$x$ is a smooth $H$-fixed point of $X$, $y = f(x)$ is a smooth 
point of $Y$, and $\Stab(y)$ is finite.  Then $i(X, x, H) = i(Y, y, H)$.
\end{thm}

\subsection*{Case I: G = H} As a first step towards proving
Theorem~\ref{thm1-1},
we will consider the case where $G = H$ is a finite abelian group.
In this case Theorem~\ref{thm1-1} can be restated as follows.

\begin{prop} \label{prop1-1.finite}
Let $H$ be a finite abelian group, 
and $f \colon X \lra Y$ be a birational morphism of 
irreducible generically free $H$-varieties of dimension $r = \rank(H)$.
Assume that $x$ is a smooth $H$-fixed point of $X$ and $y = f(x)$
is a smooth point of $Y$.  Then $i(X, x, H) = i(Y, y, H)$.
\end{prop}

Before proceeding with the proof of
Proposition~\ref{prop1-1.finite}, we introduce some background 
material on the power series ring $k[[u_1, \dots, u_r]]$.

Given $w \in k[[u_1, \dots, u_r]]$ we shall denote by $\lm(w)$
the lowest degree monomial in $u_1, \dots, u_r$ which enters into
$w(u_1, \dots, u_r)$ with a nonzero coefficient. Here ``lowest degree"
refers to a fixed (lexicographic) monomial order $\succ$ given by, say,
$u_1 \succ \dots \succ u_r$. 

Suppose $v_1, \dots, v_m$ lie in the maximal ideal 
of $k[[u_1, \dots, u_r]]$, i.e, $\lm(v_i) \succ 1$ for any $i = 1, \dots, m$. 
Then we can substitute $v_1, \dots, v_m$ into any power series 
$p \in k[[z_1, \dots, z_m]]$; in other words,
$p(v_1, \dots, v_m)$ is a well-defined element of $k[[u_1, \dots, u_r]]$.
If $p = Z$ is a monomial in $k[[z_1, \dots, z_m]]$ then clearly
\begin{equation} \label{e.mon1}
 \lm \Bigl( Z(v_1, \dots, v_m) \Bigr) = Z \Bigl( \lm(v_1), \dots, \lm(v_m)
\Bigr) \, . 
\end{equation}
We shall write $\lf< u_1, \dots, u_r \r>$ for the group of all
monomials in $u_1, \dots, u_r$ (here we allow negative exponents).

\begin{lem} \label{lem.monomial}
Suppose $v_1, \dots, v_m \in k[[u_1, \dots, u_r]]$. If
$\lm(v_1), \dots, \lm(v_m)$ generate a rank $m$ subgroup 
$\Lambda$ in $\lf< u_1, \dots, u_r \r> \simeq \bbZ^r$ then
$\lm \Bigl(p(v_1, \dots, v_m) \Bigr) \in \Lambda$ for any 
$p \in k[[z_1, \dots, z_m]]$. 
\end{lem}

Note that the conditions of the lemma imply $m \leq r$;
only the case $m = r$ will be used in the subsequent application.   

\begin{pf} Suppose
$ p(z_1, \dots, z_m) = \sum a_Z Z$, where $Z$ ranges over all monomials
in $z_1, \dots, z_m$ with non-negative exponents. 
By our assumption $\lm(v_1), \dots, \lm(v_m)$ are (multiplicatively)
linearly independent, i.e.,
\[ Z(\lm(v_1), \dots, \lm(v_m)) \neq 
 Z'(\lm(v_1), \dots, \lm(v_m)) \]
for any two distinct monomials $Z$ and $Z'$. 
Suppose $Z_0(\lm(v_1), \dots, \lm(v_m))$ has the smallest monomial order
among all monomials (in $u_1, \dots, u_m$) of the form 
$Z(\lm(v_1), \dots, \lm(v_m))$, with $c_Z \neq 0$. 
Then~\eqref{e.mon1} tells us that 
\[ \lm(Z(v_1, \dots, v_m)) \succ \lm(Z_0(v_1, \dots, v_m)) = 
Z_0 (\lm(v_1), \dots, \lm(v_m)) \] 
for any $Z \neq Z_0$ with $c_Z \neq 0$. Thus
\[ \lm(p(v_1, \dots, v_m)) = \lm(Z_0(v_1, \dots, v_m)) = 
 Z_0(\lm(v_1), \dots, \lm(v_m)) \in \Lambda \, , \]
as claimed.
\end{pf}

\begin{pf*}{Proof of Proposition~\ref{prop1-1.finite}}
Diagonalizing the action of $H$ on the cotangent space
$T^*_x(X)$, we obtain a basis $\overline{u}_1, \dots, \overline{u}_r \in
T^*_x(X)$ such that $h \overline{u}_i = \chi_i(h) \overline{u}_i$
for every $h \in H$; here $\chi_1, \dots, \chi_r \in H^*$.
Since the $k$-linear map 
$\cM_{x, X} \lra \cM_{x, X}/\cM_{x, X}^2 = T^*_x(X)$
has an $H$-invariant $k$-linear 
splitting, we can find a local system of parameters
$u_1, \dots, u_r \in \cM_{x, X}$ such that
\begin{equation} \label{e.chi}
h u_i = \chi_i(h) u_i 
\end{equation}
for every $h \in H$ and $i = 1, \dots, r$.
Similarly, we can find a local coordinate system
$v_1, \dots, v_r \in \cM_{y, Y}$ for $Y$ at $y$ and $\eta_1, \dots,
\eta_r \in H^*$ such that 
\begin{equation} \label{e.eta}
h v_i = \eta_i(h) v_i 
\end{equation}
for every $h \in H$ and $i = 1, \dots, r$.
Clearly $i(X, x, H) = \chi_1 \wedge \dots \wedge \chi_r$ and
$i(Y, y, H) = \eta_1 \wedge \dots \wedge \eta_r$.

We shall identify the elements $v_1, \dots, v_r$ with their images
in $\cO_{x, X}$ under $f^* \colon \cO_{y, Y} \hookrightarrow \cO_{x, X}$.
The $H$-action on $\cO_{x, X}$ naturally extends to
$k[[u_1, \dots, u_r]]$; in view of~\eqref{e.chi}
the leading term map $w \mapsto \lm(w)$ is $H$-equivariant. 
Suppose $\lm(v_i) = u_1^{a_{i1}} \dots u_r^{a_{ir}}$ 
for some non-negative integers $a_{ij}$. 
Then~\eqref{e.chi} and~\eqref{e.eta} imply
$\eta_i = \prod_j \chi_i^{a_{ij}}$, and thus, up to sign,
\[ i(Y, y, H) =  \eta_1 \wedge \dots \wedge \eta_r = 
\det(a_{ij}) \chi_1 \wedge \dots \wedge \chi_r = 
\pm \det(a_{ij}) i(X, x, H) \, . \]
There are two conclusions we can draw from this formula. First of all,
by Lemma~\ref{lem.stand} we know that both $i(X, x, H)$ 
and $i(Y, y, H)$ generate $\bigwedge^r(H^*)$; thus 
$\det(a_{ij}) \neq 0$.
Secondly, in order to prove the proposition, it is sufficient to show that
\begin{equation} \label{det=1}
\det(a_{ij}) = \pm 1 \; \, \text{in $\bbZ$.} 
\end{equation} 
We now proceed with the proof of~\eqref{det=1}. Let
$\lf< u_1, \dots, u_r \r>$ be the free abelian multiplicative group 
generated by $u_1, \dots, u_r$.  Since $\det(a_{ij}) \neq 0$, the leading 
monomials $\lm(v_1), \dots, \lm(v_r)$ generate a (free abelian) subgroup
$\Lambda$
of rank $r$ in $\lf< u_1, \dots, u_r \r> \simeq (\bbZ)^r$; in other words,
$\Lambda$ has finite index in $\lf< u_1, \dots, u_r \r>$. On the other hand,
\eqref{det=1} holds if and only if
$\Lambda = \lf< u_1, \dots, u_r \r>$. It is therefore
sufficient to prove that $u_i \in \Lambda$ for every $i = 1, \dots, r$.

Since $\cO_{x, X}$ and $\cO_{y, Y}$ have the same field of fractions,
each $u_i$ can be written as $p/q$, where $p, q \in \cO_{y, Y} - \{ 0 \}$.
Represent $p$ and $q$ by power 
series in $v_1, \dots, v_r$. By Lemma~\ref{lem.monomial}
$\lm(p(v_1, \dots, v_r))$ and $\lm(q(v_1, \dots, v_r))$ 
lie in $\Lambda$; thus taking the leading monomials on both sides of
the equation
\[ q(v_1, \dots, v_r) u_i = p(v_1, \dots, v_r) \, , \]
we conclude that $u_i \in \Lambda$, as desired.
\end{pf*}

\subsection*{Case II: G - arbitrary}

We are now ready to finish the proof of Theorem~\ref{thm1-1}.
The idea is to replace $X$ and $Y$ by suitable $H$-slices, then
appeal to Proposition~\ref{prop1-1.finite}. 
  
Diagonalizing the $H$-action on $T^*_y(Gy)$, we obtain a basis
$\overline{v}_1, \dots, \overline{v}_d$ such that
$h\overline v_i = \alpha_i(h)\overline v_i$ 
for some characters $\alpha_1, \dots, \alpha_d$ of $H$. Since 
the natural $H$-equivariant $k$-vector space maps
$\cM_{y, Y} \lra T^*_y(Y) \lra T^*_y(Gy)$ have $H$-equivariant splittings,
we can lift $\overline{v}_1, \dots, \overline{v}_d$ to
$v_1, \dots, v_d \in \cM_{y, Y}$ such that such that 
$h v_i = \alpha_i(h) v_i$ 
for each $h \in H$. In other words, $v_1, \dots, v_d$ form 
a local coordinate system for $Gy$ at $y$.

Since both $\Stab(x)$ and $\Stab(y)$ are finite,
$df_x \colon T_x(Gx) \lra T_y(Gy)$ is an isomorphism. Thus 
$f^*(v_1), \dots, f^*(v_r)$ form a local coordinate system for 
$Gx$ near $x$. Define $W \subset Y$ as the
irreducible component of $\{v_1 = ... = v_d = 0 \}$ passing through $y$
and $V \subset X$ as the irreducible component 
of $\{f^*(v_1) = ... = f^*(v_d) = 0 \}$ passing 
through $x$. Then
$W$ is an $H$-slice for $Y$ at $y$ and $V$ is an $H$-slice for $X$ at $x$.

Clearly $f(V) \subset W$, i.e., $f_{|V} : V \lra W$ is a well-defined 
morphism. We claim that $f_{|V}$ is, in fact, a birational morphism.
Theorem~\ref{thm1-1} follows from this claim because
\[ i(X, x, H) 
\stackrel{\text{{\tiny Remark \ref{rem2.25}}}}{=}
i(V, x, H) 
\stackrel{\text{{\tiny Proposition \ref{prop1-1.finite}}}}{=}
 i(W, y, H) 
\stackrel{\text{{\tiny Remark \ref{rem2.25}}}}{=}
 i(Y, y, H) \, . \]

To show that $f_{|V}$ is dominant, assume, to the contrary, that
$\dim(f(V)) < r$. Since $X$ is irreducible, 
$GV$ is dense in $X$ by Lemma~\ref{lem.sl3} and thus
$ \dim(Y) = \dim(f(X)) = \dim(f(GV)) = \dim(G f(V)) \le d + \dim(f(V)) 
< d + r $,
a contradiction.

It remains to show that $f_{|V}$ is generically 1-1 on closed points. 
Since $f$ is a birational morphism (i.e., has degree 1), there exists 
a dense $G$-invariant open subset $Y_0$ of $Y$ such that for every
$y _0 \in Y_0$, $f^{-1}(y_0)$ is a single point in $X$.  Since $GW$ is dense 
in $Y$ (by Lemma~\ref{lem.sl3}), $Y_0 \cap W$ is a dense open subset of $W$.
Thus a general point of $W$ has exactly one preimage in $X$.
On the other hand, a general point of $W$ has at least
$\deg(f_{|V}) \ge 1$ preimages in $X$. This shows that
$\deg(f_{|V}) = 1$, i.e. $f_{|V}$ is birational, as claimed.
\qed

\section{Proof of Theorem~\ref{thm1}}
\label{sect7}

In this section we deduce Theorem~\ref{thm1} from Theorem~\ref{thm1-1}.  
Our proof relies on canonical resolution of singularities.
(We remark that canonical resolution 
of singularities is not used elsewhere in this paper.) 

We begin with two simple preliminary results.

\begin{lem} \label{lem.Z_0} Let $H$ be an algebraic group,
$f \colon Z \lra X$ be a birational morphism of complete irreducible
$H$-varieties and $X_1$ be an irreducible
$H$-invariant codimension 1 subvariety of $X$. If $X_1$ passes through
a normal point of $X$ then there exists an $H$-invariant irreducible
codimension 1 subvariety $Z_1 \subset Z$ such 
that $f_{|Z_1} \colon Z_1 \lra X_1$ is 
a birational morphism. 
\end{lem}

\begin{pf} Since $X_1$ contains a normal point of $X$, the rational map
$f^{-1} \colon X \brokrarr Z$ is defined at a general point of $X_1$.
Now $Z_1$ = the closure of $f^{-1}(X_1)$ in $Z$, has the desired properties.
\end{pf}

\begin{prop} \label{prop.fixed} 
Let $H$ be a diagonalizable group, 
$\alpha \colon Z \lra X$ be a birational morphism of complete 
irreducible $H$-varieties and $x$ be a smooth $H$-fixed point of $X$. 
Then there exists an $H$-fixed point $z \in Z$ such that $\alpha(z) = x$.
\end{prop}

This result can be established by the argument used in the proof
of~\cite[Proposition~A.4]{ry1}, due to J. Koll\'ar and E. Szab\'o.
(In fact, if $H$ is a $p$-group then Proposition~\ref{prop.fixed} 
follows from~\cite[Proposition~A.4]{ry1}.) We give
a simple self-contained proof below.

\begin{pf} 
We argue by induction on $\dim(X)$. The base case,
$\dim(X) = 0$, is obvious.  For the induction step, assume $\dim(X) \geq 1$.
We claim that there exists codimension 1 $H$-invariant irreducible subvariety
$X_1$ such that $x$ is a smooth point of $X_1$. 
Arguing as we did at the beginning of the
proof of Proposition~\ref{prop1-1.finite}, we see that there exists 
non-zero element $u \in \cM_{x, X}$ such that for every $h \in H$,
$h u = \alpha(h) u$, where $\alpha$ is a character of $H$.
Then the (locally closed) subvariety $\{ u = 0 \}$ of $X$ is $H$-invariant 
and smooth at $x$. Hence, its unique irreducible component 
passing through $x$ is also $H$-invariant, and we can define $X_1$ 
as the closure of this irreducible component. This proves the claim.

Now by Lemma~\ref{lem.Z_0} there exists 
a codimension 1 irreducible $H$-invariant subvariety $Z_1$ of $Z$
such that $\alpha_{|Z_1} \colon Z_1 \lra X_1$ is a birational morphism
of $H$-varieties. Applying the induction assumption to this morphism,
we construct $z \in Z_1 \subset Z$ with desired properties.
\end{pf}

We are now ready to complete the proof of Theorem~\ref{thm1}.
The idea is to construct a complete smooth model $Z$
of $X$ (or $Y$) that dominates them both, i.e., fits into a diagram 
\[ \begin{array}{ccccc}
                               &   & Z & &  \\
          & \stackrel{\alpha} \swarrow   &  & \stackrel{\beta}\searrow &  \\
                           X &             & &   & Y \, . 
\end{array} \]
where $\alpha$ and $\beta$ are birational morphisms of $G$-varieties.
If we find such a $Z$, Theorem~\ref{thm1} will easily follow.
Indeed, by Proposition~\ref{prop.fixed} 
there exists an $H$-fixed point $z \in Z$ such that $\alpha(z) = x$.
Setting $y = \beta(z)$ and applying Theorem~\ref{thm1-1} to
$\alpha$ and $\beta$, we conclude that
$i(X, x, H) = i(Z, z, H) = i(Y, y, H)$, as desired.

It remains to construct $Z$.
Let $W \subset X \times Y$ be the closure of the graph of a birational
isomorphism $f$ between $X$ and $Y$. Then
$W$ is a complete $G$-variety that dominates both $X$ and $Y$. In other 
words, $W$ satisfies all of our requirements for $Z$, with one exception:
it may not be smooth. Let 
\begin{equation} \label{e.tower}
\pi \colon Z = W_n@>\pi_n>>\dots @>\pi_1>>W_0 = W\ ,
\end{equation}
be the canonical resolution of singularities of $W$,
as in~\cite[Theorem~7.6.1]{Vil} or~\cite[Theorem~1.6]{bm}. Here 
$Z$ is smooth, and each $\pi_i$ is a blowup with a smooth center;
since these centers are canonically chosen, they are $G$-invariant.
Thus the $G$ action can be lifted to $Z$ so that
$\pi$ is a birational morphism of complete $G$-varieties. 
The smooth complete $G$-variety $Z$ constructed in this way has
the desired properties.
\qed

\begin{remark} \label{rem.pf2} 
An alternative construction of $Z$ 
is given by the equivariant version of Hironaka's theorem on elimination 
of points of indeterminacy (proved in~\cite{ry4}), which asserts that
for every rational map
$f \colon X \brokrarr Y$ of $G$-varieties
there exists a sequence of blowups
\[ \pi \colon Z = X_n@>\pi_n>>\dots @>\pi_1>>X_0\ \]
with smooth $G$-equivariant centers such that 
$f \pi \colon Z \lra Y$ is regular. 
The advantage of this approach is that it only uses 
Proposition~\ref{prop.fixed} in the case where $\alpha$ is a single
blowup with a smooth $G$-equivariant center (in which case the proof 
is immediate; see, e.g.,~\cite[Lemma 5.1]{ry1}). On the other hand,
since the theorem on equivariant elimination of points of 
indeterminacy is itself deduced from canonical resolution 
of singularities in~\cite{ry4}, we opted for a direct proof here.   
\end{remark}

\begin{remark} \label{rem.irred}
A rational map $f \colon X \brokrarr Y$ (respectively, a morphism 
$f \colon X \lra Y$) of reducible varieties is called 
a birational isomorphism (respectively, a birational morphism)
if $X$ and $Y$ have irreducible decompositions
$X_1 \cup \dots \cup X_n$ and $Y_1 \cup \dots \cup Y_n$ such that
$f$ restricts to a birational isomorphism $X_i \brokrarr Y_i$ 
(respectively, a birational morphism) for each $i$.
With this definition, the irreducibility assumption in 
Proposition~\ref{prop.fixed} can be removed. Indeed, we can reduce to
the case where $X$ and $Z$ are irreducible by replacing them
with suitable irreducible components. 

The irreducibility assumption in
Theorem~\ref{thm1} (respectively, Theorem~\ref{thm1-1} and
Proposition~\ref{prop1-1.finite}) can also be removed,
if we assume $\dim_x(X) = d+r$
(respectively, $\dim_x(X) = d+ r$ and $\dim_x(X) = r$).
Indeed, if $X_1$ is the (necessarily unique) irreducible component
of $X$ containing $x$ and $G_1 = \{ g \in G \, | \, g(X_1) = X_1 \}$, 
then $H \subset G_1$, $[G:G_1]< \infty$, and $i(X, x, H) = i(X_1, x, H)$, 
so that in each case we may replace $X$, $Y$, and $G$ by $X_1$, $Y_1$ 
and $G_1$.
\end{remark} 

\begin{remark} \label{rem.thm1} One may ask if the condition that
$\Stab(x)$ is finite for every $x \in X^H$ (and similarly for $Y$)
of Theorem~\ref{thm1} is ever satisfied. Indeed, if $H$ is contained 
in a torus $T$ of $G$ then the answer is ``no", since 
$\Stab(x)$ is infinite for every $x \in X^T \subset X^H$, and 
$X^T \neq \emptyset$ by the Borel Fixed Point Theorem. On the other
hand, if the centralizer $C_G(H)$ is finite, then we claim that every
generically free $G$-variety $X$ has a birational model satisfying 
this condition. Indeed, by~\cite[Theorem 1.1]{ry1} $X$ has a model
with the property that the stabilizer of every point is of the form
$U \sdp D$, where $D$ is diagonalizable and $U$ is unipotent. 
Assume $x \in X^H$. By the Levi decomposition theorem, 
we may choose $D$ so that $D$ contains $H$. Now~\cite[Lemma 7.3]{ry1}
tells us that $U = \{ 1 \}$. Thus $\Stab(x) = D \subset C_G(H)$, which
is finite, as claimed.

Examples of pairs $H \subset G$, where $G$ is a semi-simple 
algebraic group and $H$ is an abelian subgroups of $G$ whose centralizer 
is finite can be found in~\cite{griess}; see also \cite[Section 8]{ry1}.
\end{remark}

\section{A birational classification of linear representations}
\label{sect.bir-class}

It is not difficult to see that Conjecture~\ref{conj1}
fails if $G$ is a finite cyclic group 
of order $n= 5$ or $ \geq 7$. Indeed, let
$V_{\omega}$ be the 1-dimensional
representation of $G$ such that $\sigma$ acts on $V_{\omega}$
by the character $x \mapsto \omega x$,
where $\omega$ is a primitive $n$th root of unity.  Since 
any birational automorphism of $\bbA^1$ lifts to a regular automorphism 
of $\bbP^1$, it is easy to see that $V_{\omega}$ is birationally isomorphic to 
$V_{\omega'}$ iff $\omega' = \omega$ or $\omega' = \omega^{-1}$. 
(The two $G$-fixed points in $\bbP^1$  
are preserved in the former case and interchanged in the latter.) 
If $n = 5$ or $\ge 7$, we can find two primitive $n$th roots of unity 
$\omega$ and $\omega'$ such that $\omega' \neq \omega^{\pm 1}$, so that
$V_{\omega}$ and $V_{\omega'}$ are not birationally isomorphic.
(P. I. Katsylo has informed us that this observation was 
independently made by E. A. Tevelev.) 

In this section we will classify faithful linear representations of 
diagonalizable group $G$ up to birational equivalence and
show that Conjecture~\ref{conj1} fails for a number of groups, 
both abelian and non-abelian.  The general flavor of the results 
we obtain will be similar to the situation described in the above 
paragraph but the arguments are more complicated due to the fact
that we will be working with higher-dimensional varieties, rather 
than curves. 

\subsection*{Representations of diagonalizable groups}

Recall that every linear representation $V$ of $G$ decomposes as 
a sum of 1-dimensional character spaces; if the associated characters
of $G$ are $\chi_1, \dots, \chi_d$, we shall write $V = \chi_1 \oplus
\dots \oplus \chi_d$. 

\begin{thm} \label{thm6.1} Let $G$ be a diagonalizable group of
rank $r$ and let $V = \chi_1 \oplus \dots \oplus \chi_d$ 
and $W = \eta_1 \oplus \dots \oplus \eta_d$ be faithful $d$-dimensional
linear representations of $G$. (In particular, $d \geq r$.) Then 
$V$ and $W$ are birationally isomorphic as $G$-varieties if and only if
$\chi_1 \wedge ... \wedge \chi_d =  \pm \eta_1 \wedge \dots \wedge
\eta_d$ in $\bigwedge^d(G^*)$.
\end{thm}

\begin{pf} 
Since $G$ acts faithfully on $V$ and $W$, we have 
\begin{equation} \label{e6.1}
\lf< \chi_1, \dots, \chi_d \r> = 
\lf< \eta_1, \dots, \eta_d \r> = G^*  \, .
\end{equation}

Assume $\chi_1 \wedge ... \wedge \chi_r =  \pm \eta_1 \wedge \dots \wedge
\eta_r$. Then by Corollary~\ref{cor.equiv} there exists an $N = (n_{ij})
\in \GL(\bbZ)$ such that $\eta_i = \chi_1^{n_{1i}} \dots \chi_d^{n_{di}}$.
The desired birational isomorphism $V \brokrarr W$ is can now be
explicitly given, in the
natural coordinates on $V$ and $W$, by
$(x_1, \dots, x_d) \lra (y_1, \dots, y_d)$, where
$y_i = x_1^{n_{1i}} \dots x_d^{n_{di}}$.

Conversely, suppose
\begin{equation} \label{e.chi-neq-eta} \textstyle
\chi_1 \wedge ... \wedge \chi_d \neq  \pm \eta_1 \wedge \dots \wedge
\eta_d \; \, \text{in} \; \, \bigwedge^d (G^*)
\end{equation}
We want to prove that $V$ and $W$ are not birationally isomorphic
as $G$-varieties.  
Note that~\eqref{e.chi-neq-eta} is impossible if $d \geq r+1$, since 
in this case
$\bigwedge^d (G^*) = (0)$; see Lemma~\ref{lem.ab1}(b). Thus we will assume from
now on that $d = r = \rank(G)$.  We will consider three cases.

\smallskip
Case 1: $G = (G_m)^r$ is a torus. In this 
case $\chi_1 \wedge ... \wedge \chi_r$ and
$\eta_1 \wedge \dots \wedge \eta_r$ are both generators of
$\bigwedge^r (G^*) = \bbZ$, so that~\eqref{e.chi-neq-eta} is impossible.

\smallskip
Case 2: $G$ is a finite abelian group.  
The $G$-action on $V = \bbA^r$ (respectively $W = \bbA^r$)
naturally extends to the projective space $\overline{V} = \bbP^r$
(respectively, $\overline{W} = \bbP^r$).
Example~\ref{ex.p^r} shows that for every $G$-fixed point $x \in \overline{V}$, 
$i(\overline{V}, x, G) = \pm \chi_1 \wedge \dots \wedge \chi_r$ and for every
$G$-fixed point $y \in \overline{W}$, 
$i(\overline{W}, y, G) = \pm \eta_1 \wedge \dots \wedge \eta_r$.
Thus in view of~\eqref{e.chi-neq-eta}, Theorem~\ref{thm1}
says that $\overline{V}$ and $\overline{W}$ (and, hence, $V$ and $W$)
are not birationally isomorphic as $G$-varieties. 

\smallskip
Case 3: $G$ is a diagonalizable group but not a torus. 
Write $G = \bbG_m(n_1) \times \dots \times \bbG_m(n_r)$, as
in~\eqref{e.fg-abel}, with $n_1 \geq 2$.
Let $H = \bbG_m(n_1)^r = (\bbZ/n_1\bbZ)^r$ be the $n_1$-torsion subgroup of
$G$.  It is sufficient to show that 
$V$ and $W$ are not birationally isomorphic as $H$-varieties; then
they certainly cannot be birationally isomorphic as $G$-varieties. 
By Case 2, it is enough to show that
\begin{equation} \label{e.primed-char} \textstyle
\chi_1' \wedge ... \wedge \chi_r' \neq  \pm \eta_1' \wedge \dots \wedge
\eta_r' \; \, \text{in} \; \, \bigwedge ^r(H^*)
\end{equation}
where $\chi_i'$ and $\eta_j'$ are the characters of $H$ obtained by 
restricting $\chi_i$ and $\eta_j$ from $G$ to $H$. Note that
the inclusion $\phi \colon H \hookrightarrow G$ induces a surjection
$\phi^* \colon G^* \lra H^*$ of the dual group, which, in turn,
induced a map of cyclic groups
$\bigwedge^r(\phi^*) \colon \bigwedge^r(G^*) \lra \bigwedge^r(H^*)$.
Elementary group theory tells us that
$G^* = (\bbZ/n_1\bbZ) \times \dots \times (\bbZ/n_r \bbZ) $,
$H^* = (\bbZ/n_1\bbZ)^r$, $\phi^* \colon G^* \lra H^*$ is 
(componentwise) reduction modulo $n_1$, and
$\bigwedge^r(\phi^*)$ is the identity map
$\bigwedge^r(G^*) = \bbZ/n_1\bbZ \stackrel{\simeq}{\lra} \bbZ/n_1\bbZ
= \bigwedge^r(H^*)$. Applying 
$\bigwedge^r(\phi^*)$ to both sides of~\eqref{e.chi-neq-eta}, we 
obtain~\eqref{e.primed-char}, as desired. 
\end{pf}

\subsection*{Proof of Theorem~\ref{thm6.2}}
Let $G$ be a diagonalizable group $G$ of rank $r$ of the 
form~\eqref{e.fin-abel}, and let
$V = \chi_1 \oplus \dots \oplus \chi_d$ 
and $W = \eta_1 \oplus \dots \oplus \eta_d$ be faithful $d$-dimensional
linear representations of $G$.

\smallskip
(a) If $d \geq r + 1$ then $\bigwedge^r(G^*) = (0)$, so that
$\chi_1 \wedge ... \wedge \chi_d =  0 = \eta_1 \wedge \dots \wedge
\eta_d$.  Thus
$V$ and $W$ are birationally isomorphic as $G$-varieties by
Theorem~\ref{thm6.1}.

\smallskip
 From now on we will assume $d = r$.
Note that in this case both $\chi_1 \wedge ... \wedge \chi_d$ 
and $\eta_1 \wedge \dots \wedge \eta_d$ are generators of
$\bigwedge^r(G^*) = \bbZ/n_1\bbZ$. 

\smallskip
(b) Suppose $n_1 = 2$. Since $\bigwedge^r(G^*) = \bbZ/2\bbZ$
has only one generator, 
$\chi_1 \wedge ... \wedge \chi_d =  \eta_1 \wedge \dots \wedge \eta_d$. 
Thus $V$ and $W$ are birationally isomorphic as $G$-varieties by
Theorem~\ref{thm6.1}.

Now assume $n_1 = 0$.  Then $\bigwedge^r(G^*) = \bbZ$. The only generators
of this group are $\pm 1$; thus
$\chi_1 \wedge ... \wedge \chi_d = \pm  \eta_1 \wedge \dots \wedge \eta_d$, 
and, once again, Theorem~\ref{thm6.1} tells us that
$V$ and $W$ are birationally isomorphic. 

\smallskip
(c) Suppose $n_1 \geq 3$.
By Theorem~\ref{thm6.1}, birational isomorphism classes of 
$r$-dimensional linear representations of $H$ are in 1-1 correspondence 
generators of $\bigwedge^r(H^*) \simeq \bbZ/n_1\bbZ$ 
(as an additive group), modulo multiplication by $-1$. 
Since $n_r \geq 3$, $a \neq -a$ for any
generator $a$ of $\bbZ/n_1\bbZ$.  Thus in this case  the number 
of isomorphism classes of faithful $r$-dimensional $H$-representations 
is $\phi(n_r)/2$, as claimed.
\qed

\subsection*{Further counterexamples to Conjecture~\ref{conj1}}

Theorem~\ref{thm6.2} shows that Conjecture~\ref{conj1} fails 
for many diagonalizable groups. We will now see this conjecture
fails for some non-abelian groups as well.

\begin{prop} \label{prop8.1}
Let $n$ and $r$ be a positive integers, $P$ be a subgroup
of $\Sym_r$ and $G = (\bbZ/n\bbZ)^r \sdp P$, where $P$ acts on 
$(\bbZ/n\bbZ)^r$ by permuting the factors.  Assume there 
exists an $m \in \bbZ$ such that $(m, n) = 1$ and 
$m^r \not \equiv \pm 1 \pmod n$.  Then there exist two
birationally inequivalent $r$-dimensional representations of $G$. 
In particular, Conjecture~\ref{conj1} fails for this group.
\end{prop}

We remark that an integer $m$ satisfying the requirements of 
Proposition~\ref{prop8.1} always exists if the exponent 
of $U_n$ does not divide $2r$; here $U_n$ is the (multiplicative) group 
of units in $\bbZ/n\bbZ$. In particular, $m$ exists  
if there is a prime power $p^e$ such that $p^e \, | \, n$ but
$\phi(p^e) = (p-1)p^{e-1} \notdiv 2r$.

\begin{pf} Let $\omega$ be a primitive $n$th root of unity.
We define the $r$-dimensional representations $V$ and $W$ of $G$
as follows:
\begin{alignat*}{2}
 ((a_1, \dots, a_r), \sigma) \colon && (v_1, \dots, v_r) & \mapsto
(\omega^{a_1} v_{\sigma^{-1}(1)}, \dots,
\omega^{a_r} v_{\sigma^{-1}(r)}) \\
\intertext{and}
 ((a_1, \dots, a_r), \sigma) \colon && (w_1, \dots, w_r) & \mapsto
(\omega^{ma_1} w_{\sigma^{-1}(1)}, \dots,
\omega^{ma_r} w_{\sigma^{-1}(r)}) \, . 
\end{alignat*}
Here $a_1, \dots, a_n \in \bbZ/n\bbZ$, $\sigma \in P \subset \Sym_r$, 
$(v_1, \dots, v_r) \in V$ and $(w_1, \dots, w_r) \in W$.
It is easy to see that $V$ and $W$ are, indeed, well-defined 
faithful $r$-dimensional representations of $G$. 

To prove 
the proposition it is sufficient to show that $V$ and $W$ are not
birationally isomorphic as $(\bbZ/n\bbZ)^r$-varieties.
Let $\chi_i$ be the character of $(\bbZ/n\bbZ)^r$ given by 
$\chi(a_1, \dots, a_r) = \omega^{a_i}$. Then, as 
$\bbZ/n \bbZ$-representations, 
$V = \chi_1 \oplus \dots \oplus \chi_r$ and
$W = \chi_1^m \oplus \dots \oplus \chi_r^m$. 
By our assumption
\[ \chi_1^m \wedge \dots \wedge \chi_r^m = 
m^r \chi_1 \wedge \dots \wedge \chi_r \neq
\pm \chi_1 \wedge \dots \wedge \chi_r \]
in $\bigwedge^r(\bbZ/n\bbZ^*) \simeq \bbZ/n\bbZ$.
Thus Theorem~\ref{thm6.1} tells us that $V$ and $W$ are not 
isomorphic as $(\bbZ/n\bbZ)^r$-varieties (and hence, as $G$-varieties).
\end{pf}

\begin{remark} \label{rem8.2} 
The same argument proves the following stronger result.
Let $n_1, n_2, \dots, n_s$, $r_1, \dots, r_s$ be positive
integers such that $n_i$ divides $n_{i+1}$ for $i = 2, \dots, r$,
let $P_i$ be a subgroup of the symmetric group $\Sym_{r_i}$ and let
$G_i = (\bbZ/n_i \bbZ)^{r_i} \sdp P_i$. Assume there exist integers
$m_1, \dots, m_s$ such that $(m_i, n_i) = 1$ and 
$m_1^{r_1} \dots m_s^{r_s} \not \equiv \pm 1 \pmod{n_1}$. 
Then $G = (G_m)^a \times G_1 \times \dots \times G_s$
has two birationally inequivalent
($a + r_1 + \dots + r_s$)-dimensional representations.
In particular, Conjecture~\ref{conj1} fails for any $G$ of this form.
\end{remark}

\begin{remark} \label{rem8.3} 
The proof of Proposition~\ref{prop8.1} shows that
$G = (\bbZ/n\bbZ)^r \sdp P$ has at least $|\pm U_n^r|/2$ 
birational isomorphism classes of $r$-dimensional 
representations.  Here, as before, $U_n$
denotes the multiplicative group of units in the ring $\bbZ/n\bbZ$, 
and $\pm U_n^r$ denotes the subset of $U_r$ consisting of elements of
the form $\pm m^r$, as $m$ ranges over $U_n$.

A similar estimate can be given for the number of birational isomorphism 
classes of ($a + r_1 + \dots + r_s$)-dimensional representations of the 
group $G$ in Remark~\ref{rem8.2}. In particular, if $n_1 = \dots = n_s$
and $(r_1, \dots, r_s) = 1$ then there are at least $\phi(n_1)/2$ such 
classes.
\end{remark}

\section{Birational equivalence of quantum tori}
\label{sect.pgl_n}

In this section we will use the invariant $i(X, x, H)$ to classify
$\PGL_n$-varieties (and consequently central simple algebras) of a certain 
form. In particular, we will prove Theorem~\ref{thm3}.

\subsection*{Abelian subgroups of $\protect \PGLn$}

Let $A$ be a finite abelian group of order $n$ and let $V = k[A]$
be the group ring of $A$. For $a \in A$ and $\chi \in A^*$
define $P_a, D_{\chi} \in \GL(V)$ by $P_a(b) = ab$ and
$D_{\chi}(b) = \chi(b)b$ for every $b \in A$. It is easy to see that
$D_{\chi}P_a = \chi(a) P_a D_{\chi}$. Thus if $p_a$ and $d_{\chi}$ denote
the elements of $\PGL(V)$ represented, respectively, by 
$P_a$ and $D_{\chi} \in \GL(V)$, then 
\begin{align*} \label{e.phi}
\phi \colon A \times A^* & \hookrightarrow \PGL(V) = \PGLn \\
(a, \chi) & \lra p_ad_{\chi}
\end{align*}
defines an embedding of $A \times A^*$ in $\PGLn$.

Let $H$ be an abelian subgroup of $\PGLn$. Then $H$ is naturally
equipped with an alternating bilinear form $\omega_H \colon
H \times H \lra \mu_n$ (cf. Definition~\ref{def.sympl}(a)). 
Here
$\mu_n$ is the group of $n$th roots of unity in $k$, identified with the
center of $\SL_n(k)$, and
$\omega_H(a, b) = ABA^{-1}B^{-1}$, where
$A$ and $B \in \SL_n$ represent $a$ and $b \in \PGLn$ respectively.

\begin{lem} \label{lem.sympl-H}
Let $A$ be a finite abelian group of rank $r$,
$H = \phi(A \times A^*) = \{p_a d_{\chi} \, | \, a \in A , \; \chi \in A^* \}$ 
be the subgroup of $\PGLn$ defined above. Then 

\smallskip
(a) the elements of $P_a D_{\chi}$ span $\Mn$ as a $k$-vector space,
as $a$ ranges over $A$ and $\chi$ ranges over $A^*$, and

\smallskip
(b) the alternating bilinear form $\omega_H$ is symplectic 
(i.e., non-degenerate).

\smallskip
\noindent
Let $g$ be an element of the normalizer 
$N_{\PGLn}(H)$, and $\psi_g \colon H \lra H$ be conjugation by $g$. 
Then

\smallskip
(c) $\psi_g$ preserves $\omega_H$, and  

\smallskip
(d) $\psi_g$ induces the identity map 
$\bigwedge^{2r}(H^*) \lra \bigwedge^{2r}(H^*)$.
\end{lem}

\begin{pf} 
(a) See~\cite[Lemma 3.2]{tsch}.  $\quad$ 
(b) See~\cite[Lemma 7.8]{ry2}. 

\smallskip
(c) Choose $a$ and $b \in H \subset \PGLn$ and lift them to
$A$ and $B \in \SL_n$. Since $ABA^{-1}B^{-1}$ is a central
element of $\SL_n$, we have 
\begin{align*}
\omega_H(\psi_g(a), \psi_g(b)) & = \omega_H (gag^{-1}, gbg^{-1}) \\
& = (gAg^{-1})(gBg^{-1})(gA^{-1}g^{-1}) (gB^{-1}g^{-1}) \\
& = g (ABA^{-1}B^{-1}) g^{-1} = ABA^{-1}B^{-1} = \omega_H(a, b) \, , 
\end{align*}
as claimed.

\smallskip
(d) Follows from (b), (c) and Lemma~\ref{lem.sf2}(b).
\end{pf}

\subsection*{$\protect \PGLn$-varieties}

\begin{prop} \label{prop.pgl_n} Let $A$ be a finite abelian group of
order $n$ and rank $r$ and let $H = \phi(A \times A^*)$ be the subgroup
of $\PGLn$ defined above.
Suppose $V = \chi_1 \oplus \dots \oplus \chi_{2r}$ and
$W = \eta_1 \oplus \dots \oplus \eta_{2r}$ are faithful
representations of $H$. Then the following are equivalent:

\smallskip
(a)  $\chi_1 \wedge \dots \wedge \chi_{2r} 
 = \pm \eta_1 \wedge \dots \wedge \eta_{2r}$ in $\bigwedge^{2r}(H^*)$,

\smallskip
(b) $V$ and $W$ are birationally isomorphic as $H$-varieties, 

\smallskip
(c) $X = \PGL_n *_H V$ and $Y = \PGLn *_H W$ are birationally 
isomorphic as $\PGLn$-varieties. 
\end{prop}

\noindent
Here $\PGL_n *_H V$ and $\PGLn *_H W$ are homogeneous fiber products; 
see Example~\ref{ex.induced}.

\begin{pf} (a) and (b) are equivalent by Theorem~\ref{thm6.1}. 
The implication (b) $\Longrightarrow$ (c) is obvious.

Thus we only need to show (c) $\Longrightarrow$ (a). The idea of 
the proof is to appeal to Theorem~\ref{thm1}. We begin by 
observing that $X$ and $Y$ naturally embed as dense 
open subsets in projective varieties $\overline{X} = 
(\bbP(\Mn) \times \overline{V})/H$ and $\overline{Y} = 
(\bbP(\Mn) \times \overline{W})/H$ respectively.
Here $\overline{V} = \bbP^{2r}$ 
is the projective completion of $V = \bbA^{2r}$; $\PGLn$ acts
on $\bbP(\Mn) \times \overline{V}$ by left multiplication on the first
factor; this action commutes with the $H$-action on  
$\bbP(\Mn) \times \overline{V}$ given by 
$h \colon (x, y) \lra (xh^{-1}, hy)$ and thus descends to the geometric
quotient $\overline{X} =  (\bbP(\Mn) \times \overline{V})/H$. 
We shall denote the point $x \in \overline{X}$ corresponding to 
the orbit of $(g, v) \in \bbP(\Mn) \times \overline{V}$ by
$[g, v]$.  The $H$-variety $\overline{W}$ 
and the $\PGL_n$-variety $\overline{Y}$ are defined in 
a similar manner.

Our goal is to show that

\smallskip
(i) every $H$-fixed points of $\overline{X}$ is of the form
$x = [g, v]$, where $g \in N_{\PGLn}(H)$ and $v$ is an $H$-fixed point of 
$\overline{V}$, and for any such point $x$,

\smallskip
(ii) $\Stab(x) = H$ and

\smallskip
(iii) up to sign, $i(X, x, H) = \chi_1 \wedge \dots \wedge \chi_{2r}$. 

\smallskip
\noindent
These assertions, in combination with Theorem~\ref{thm1}, 
will prove that if $X$ and $Y$ (and hence,
$\overline{X}$ and $\overline{Y}$) are birationally isomorphic then
$\chi_1 \wedge \dots \wedge \chi_{2r} = 
\pm \eta_1 \wedge \dots \wedge \eta_{2r}$, i.e.,
(c) $\Longrightarrow$ (a). 

\smallskip
To prove (i), assume $x = [g, v]$ is an $H$-fixed point of $\overline{X}$
for some $g \in \bbP(\Mn)$ and $v \in \overline{V}$.
This means that for every $h \in H$ there exists an $h' \in H$ such that
$(hg, v) = (gh', (h')^{-1}v)$ in $\bbP(\Mn) \times \overline{V}$.
Equivalently, $hg = g h'$ and $(h')^{-1}v = v$.

\newcommand{\RKer}{\operatorname{RKer}}
Consider the vector space $k^n$ of $(n \times 1)$-row vectors.
The multiplication by $g$ on the right yields a linear map
$k^n \lra k^n$; let $\RKer(g)$ be the kernel of this map.
Note that since
$g \in \bbP(\Mn)$, this linear map is only defined up to a nonzero
constant multiple but $\RKer(g)$ is well-defined.

The equality $hg = gh'$ implies that $\RKer(g)$ is 
an $H$-invariant subspace of $k^n$ with respect to the right action of
$H$; again, as $H\subset\PGLn$, the right multipication by an element
$h\in H$ is a linear map $k^n \lra k^n$ defined up to a nonzero
constant multiple, but the notion of $H$-invariance of a linear
subspace of $k^n$ is well-defined.

Now recall that by Lemma~\ref{lem.sympl-H}(a)
the $n^2$ elements of the form $P_a D_{\chi}$ which represent
the elements of $H \subset \PGL(V) = \PGLn$ in $\GL(V) = \GL_n$ span
$\Mn$ as a $k$-vector space.  Thus the only $H$-invariant subspaces 
of $k^n$ are the ones that are invariant under all of $\Mn$, namely
$k^n$ and $(0)$.  If $\RKer(g) = k^n$
then $g$ is the zero matrix, which is impossible for any
$g \in \bbP(\Mn)$. Thus we conclude that
$\RKer(g) = (0)$. This means that $g$ is nonsingular, i.e., 
$g \in \PGLn$. Now we can rewrite $hg = gh'$ as
$g^{-1}hg = h' \in H$; this shows that $g \in N_{\PGLn}(H)$. 
Moreover, as $h$ ranges over $H$, $h' = g^{-1}hg$ also 
ranges over all of $H$. Thus the equality $(h')^{-1}v= v$ implies that $v$ 
is an $H$-fixed point of $\overline{V}$. This proves (i).
 
\smallskip
\noindent
 From now on let $x = [g, v]$ be an $H$-fixed point of $\overline{X}$,
where $g \in N_{\PGLn}(H)$ and $v$ is an $H$-fixed point 
of $\overline{V}$. 

\smallskip
To prove (ii), assume  $g' \in  \Stab(x)$, i.e., $g'[g, v] = [g, v]$.
Then $g'g = g h'$ for some $h' \in H'$. Since $g \in N_{\PGLn}(H)$,
we conclude that $g' = ghg^{-1} \in H$, as desired.

\smallskip
To prove (iii), first note that $i(\overline{X}, [g, v], H) = \wedge^{2r}
\psi_g^* \Bigl( i(\overline{X}, [1, v], H) \Bigr)$, where
$\psi_g \colon H \lra H$ is conjugation by $g \in N_{\PGLn}(H)$ and
$\wedge^{2r} (\psi_g^*)$ is the automorphism of 
$\bigwedge^{2r}(H^*)$ induced by $\psi_g$; 
see Remark~\ref{rem2.4}.
By Lemma~\ref{lem.sympl-H}(d) 
$\wedge^{2r} \psi_g^*$ is the identity automorphism.  Thus
$i(\overline{X}, [g, v], H) = i(\overline{X}, [1, v], H)$. 
On the other hand, by Remark~\ref{rem2.25}
$i(\overline{X}, [1, v], H) = i(\overline{V}, v, H)$.
Finally, recall that for any $v \in \overline{V}^H$,
$i(\overline{V}, v, H) = i(V, 0_V, H) = 
\chi_1 \wedge \dots \wedge \chi_{2r}$; 
see Example~\ref{ex.p^r}. In summary,
\[ i(\overline{X}, [g, v], H) = i(\overline{X}, [1, v], H) = 
i(\overline{V}, v, H) = 
\chi_1 \wedge \dots \wedge \chi_{2r} \, , \]
as claimed.
\end{pf}

\begin{remark} \label{rem.E_8}
Recall that the exceptional group $E_8$ has a unique 
nontoral subgroup isomorphic to $(\bbZ/5\bbZ)^3$ (up to conjugation).
Denote this subgroup by $H$. Then, modifying the proof of 
Proposition~\ref{prop.pgl_n}, we can show the following:

\smallskip
{\em Let $V = \chi_1 \oplus \chi_2 \oplus \chi_3$ and 
$W = \eta_1 \oplus \eta_2 \oplus \eta_3$ be faithful 3-dimensional
representations of $H$, where $\chi_i$ and $\eta_j$ are characters
of $H$. Then the following are equivalent:

\smallskip
(a) $\chi_1 \wedge \chi_2 \wedge \chi_3 = \pm \eta_1 \wedge \eta_2 
\wedge \eta_3$ in $\bigwedge^3(H^*) \simeq \bbZ/5\bbZ$,

\smallskip
(b) $V$ and $W$ are birationally isomorphic as $H$-varieties, and

\smallskip
(c) $E_8 *_H V$ and $E_8 *_H W$ are birationally
isomorphic as $E_8$-varieties.}

\smallskip
\noindent
In particular, there are exactly two birational
isomorphism classes of $E_8$-varieties
of the form $E_8 *_H V$, where $V$ is a faithful 3-dimensional
representation of $H$: one corresponds to $\pm 1$, and the other
to $\pm 2$ in $\bbZ/5\bbZ \simeq \bigwedge^3(H^*)$.  
\end{remark}

\begin{remark} Note that the $\PGLn$-varieties $X = \PGL_n *_H V$ 
and $Y = \PGLn *_H W$ of Proposition~\ref{prop.pgl_n}, are stably
isomorphic. In fact, $X \times \bbA^1 \simeq Y \times \bbA^1$ 
because $X \times \bbA^1 = \PGL_n *_H (V \times \bbA^1)$,
$Y \times \bbA^1 = \PGL_n *_H (W \times \bbA^1)$, and $V \times \bbA^1 \simeq 
W \times \bbA^1$ as $H$-varieties by Theorem~\ref{thm6.1}.  

For the same reason $X \times \bbA^1$ and $Y \times \bbA^1$ 
are isomorphic $E_8$-varieties, if $X$ and $Y$ are 
as in Remark~\ref{rem.E_8}.
\end{remark}

\subsection*{Proof of Theorem~\ref{thm3}}

Recall that birational isomorphism classes of generically free
irreducible $\PGLn$-varieties $X$ with $k(X)^{\PGLn} = K$ are in
1-1 correspondence with central simple algebras 
of degree $n$ over $K$; see~e.g.,~\cite[X.5]{serre}
or~\cite[Section 3]{ry2}. 

In particular, by~\cite[Lemma 4.2]{tsch}, the algebra 
$Q(\omega_1, \dots, \omega_r)$ of Theorem~\ref{thm3} corresponds 
to the variety $X = \PGLn *_H V$, where $V$ is a faithful 
$2r$-dimensional representations of  $H = A \times A^*$
constructed as follows. 
Choose a set of generators $a_1, \dots, a_r$ 
for $A = \bbZ/n_1\bbZ \times \dots \times \bbZ/n_r \bbZ$ and 
a ``dual" set of generators $\chi_1, \dots, \chi_r$ for $A^*$ so that 
\[ \chi_i(a_j) = \left\{ \begin{array}{ll} 1 & \text{\rm if $i \neq j$} \\
                         \omega_i & \text{\rm if $i = j$} \, , 
\end{array} \right. \]
Now note that each $a \in A$ defines a character of $H = A \times A^*$  
by $(b, \eta) \mapsto \eta(a)$. Similarly, each
$\chi \in A^*$ gives rise to a character $H = A \times A^* \lra k^*$ 
via $(b, \eta) \mapsto \chi(b)$; we shall denote these characters
by $c(a)$ and $c(\chi)$ respectively. In these notations,
\[ V = c(a_1) \oplus \dots \oplus c(a_r) \oplus c(\chi_1)^{-1} \oplus \dots 
\oplus c(\chi_r)^{-1} \, . \]
see~\cite[Proof of Lemma 4.2]{tsch}.

Similarly the $\PGLn$-variety associated to
$Q(\omega_1^{m_1}, \dots, \omega_r^{m_r})$ is
$Y = \PGLn *_H W$, where 
$W = c(a_1') \oplus \dots \oplus c(a_r') \oplus c(\chi_1')^{-1} 
\oplus \dots \oplus c(\chi_r')^{-1}$. Here
$a_1', \dots, a_r'$ are generators of $A$ 
and $\chi_1', \dots, \chi_r'$ are generators of $A^*$ such that
\[ \chi_i'(a_j') = \left\{ \begin{array}{ll} 1 & \text{\rm if $i \neq j$} \\
                         \omega_i^{m_i} & \text{\rm if $i = j$} \, . 
\end{array} \right. \]
A natural choice for $a_i'$ and $\chi_i'$ is
$a_i' = a_i$ and $\chi_i' = \chi_i^{m_i}$, so that
\[ W = c(a_1) \oplus \dots \oplus c(a_r) \oplus c(\chi_1)^{-m_1} \oplus \dots 
\oplus c(\chi_r)^{-m_r} \, . \]

As we mentioned above, $Q(\omega_1, \dots, \omega_r)$ and 
$Q(\omega_1^{m_1}, \dots, \omega_r^{m_r})$ are isomorphic as $k$-algebras
iff their associated $\PGLn$-varieties,
$X = \PGL_n *_H V$ and $Y = \PGLn *_H W$, are birationally isomorphic. 
By Proposition~\ref{prop.pgl_n} $X$
and $Y$ are birationally isomorphic iff 
\begin{eqnarray*} 
\lefteqn{c(a_1) \wedge \dots \wedge c(a_r) \wedge c(\chi_1)^{-1} \wedge \dots
\wedge c(\chi_r)^{-1} =} \\
 & & \pm
c(a_1) \wedge \dots \wedge c(a_r) \wedge c(\chi_1)^{-m_1} \wedge \dots
\wedge c(\chi_r)^{-m_r} \; \, 
\text{in $\textstyle \bigwedge^{2r}(H^*) \simeq \bbZ/n_r \bbZ$.} 
\end{eqnarray*}
The last condition
is equivalent to $m_1 \dots m_r = \pm 1 \pmod{n_1}$.
\qed

\end{document}